\documentclass[12pt]{amsart}
\usepackage[headings]{fullpage}
\usepackage{longtable}
\usepackage{rotating}
\usepackage{amssymb}
\usepackage{xcolor}
\usepackage{hyperref}

\theoremstyle{plain}
\newtheorem{theorem}{Theorem}[section]

\theoremstyle{remark}

\theoremstyle{definition}

\newcommand{\PP}{\mathbb P}
\newcommand{\lra}{\longrightarrow}
\newcommand\CC{{\mathbb C}}

\numberwithin{table}{section}
\begin{document}
	\leavevmode
	\def\wl#1{\hbox to 10mm{\hfill$\ifcase#1\or I_{2}\or I_{4}\or I_{0}^{*}\or I_{2}^{*}\fi$\hfill}}
	

	\title[Elliptic and K3 fibrations on double octic Calabi-Yau threefolds]
	{Elliptic and K3 fibrations on double octic Calabi-Yau threefolds}
	\author{S{\l}awomir Cynk}
	\address{Institute of Mathematics, Jagiellonian University,
		{\L}ojasiewicza 6, 30-348 Krak\'ow, Poland 
	}
	
	\email{slawomir.cynk@uj.edu.pl}
	\author{Beata Kocel--Cynk}
	\address{Department of Applied Mathematics,
		Faculty of Computer Science and Telecommunications, 
		Cracow University of Technology,
		ul.~Warszawska~24,
		31-155~Krak\'ow,
		POLAND}
	\email{beata.kocel-cynk@pk.edu.pl}
	\thanks{Authors were partially supported by National Science Centre, Poland, grant No. 2020/39/B/ST1/03358.}
	
	\keywords{Calabi--Yau threefold, double octic, automorphism, K3 fibration}
	\subjclass[2010]{Primary: 14J32, Secondary: 14D06}
	
	\maketitle
	\begin{abstract}
		We study fibrations by elliptic curves and K3 surfaces of double octic Calabi-Yau threefolds determined by singular lines and points of multiplicity at least 4 of the defining octic arrangement. 
		As a consequence we conclude that every double octic admits many elliptic and K3 fibrations.
		We apply descriptions of some K3 fibrations to construct a birational map between elements of two families of double octics.
		
	\end{abstract}
	
	\section*{Introduction}
	In this note we study fibrations on double octic Calabi-Yau threefolds, by a Calabi-Yau threefold we mean a smooth, complex, projective manifold $X$ with trivial canonical line bundle $\omega_{X}\cong\mathcal O_{X}$ and vanishing first cohomology of the structure sheaf $H^{1}(\mathcal O_{X})=0$.
	A double octic is a Calabi-Yau threefold which is a resolution of singularities of a double covering of the projective space $\PP^{3}$ branched along a union of eight planes (satisfying some mild assumptions). 
	
	Double octics form a relatively small but very suitable for explicit studies family of Calabi-Yau threefolds. 
	They provide a rich source of examples of Calabi-Yau threefilds with special properties including: Calabi-Yau threefold in characterictic 3 with no lifting to characteristic zero, counter-example to Bogomolov-Tian-Todorov in positive characteristic (\cite{CvS3}),
	example of Calabi-Yau threefold with trivial monodromy and  no smooth filling (\cite{CvS5}), counter example to N\'eron-Ogg-Shafarevich criterion for Calabi-Yau threefolds (\cite{C-O}), examples of one parameter families of Calabi-Yau threefolds without a point of Maximal Unipotent Monodromy (\cite{CvS4}), examples of Hilbert modular Calabi-Yau threefolds (\cite{CScvS}). In several of these construction a crucial role plays an explicit description of an elliptic or Kummer fibrations.

	If $X\lra B$ is a fibration of a Calabi-Yau threefold, then the generic fiber $X_{b}:=F^{-1}(b)$	is of one of the following three types: elliptic curve, K3 surface abelian surface. P.M.H.~Wilson gave a characterization of the three possibilities in terms of the k\"ahler cone (\cite{wilson2}) and a condition for existence of an elliptic fibration on a Calabi-Yau threefold (\cite{wilson2}). Oguiso and Peternell (\cite{OP}) proved that on a Calabi-Yau  threefold there exists only finite number of elliptic or K3 fibrations. 
	
	Existence of elliptic or K3 fibrations for a Calabi-Yau manifold have singificant consequences for its geometry (\cite{GW, HK, KMW}), elliptic Callabi-Yau threefolds play special role in applications in physics (\cite{LY,MT}). Majority of examples originating from best known constructions  of Calabi-Yau threefolds admits an elliptic fibration (f.i. see \cite{HT}).
	In fact it is expected that any Calabi-Yau threefold of sufficiently large Picard number admits an elliptic fibration.

	In the section 2. we briefly recall properties of double octic Calabi-Yau threefolds. In the further sections we describe the following type of fibrations on a double octic
	\begin{itemize}
		\item an elliptic fibration determined by lines through an arrangement point $P$ of multiplicity at least 4,
		\item an elliptic fibration determined by lines intersecting two disjoint double or triple arrangement lines,
		\item a Kummer fibration determined by a pair of fourfold points,
		\item a K3 fibration determined by planes through a double or triple arrangement line,
		\item a K3 fibration determined by a pencil of quadrics containing two pairs of double or triple arrangement line.
	\end{itemize} 
	
	As an application of fibrations of the last type we find a birational map between elements of families of double octics No. 35 and No. 71.	
	
	\section{Double octic Calabi-Yau manifold}
	A double octic is a Calabi-Yau threefold defined as a (crepant) resolution of singularities of a double covering of the projective space $\PP^{3}$ branched along a union of eight planes 
	\[D=P_{1}\cup\cdots\cup P_{8}.\]
	If we denote by $F_{i}(x,y,z,t)$ the linear equation of the plane $P_{i}$ then the singular double octic is given in the weighted projective space $\PP(1,1,1,1,4)$ by the following equation 
	\[u^{2}=F_{1}(x,y,z,t)\cdots F_{8}(x,y,z,t).\]
	If the eight planes constituting an octic arrangement $D$ satisfies the following two conditions
	\begin{itemize}
		\item [] no six intersect in a point, 
		\item [] no four intersect along a line,
	\end{itemize}
	then this double cover admits a resolution of singularities that is a Calabi-Yau threefold (which we call a double octic). An explicit procedure of a resolution of singularities is described in \cite{CSz}. 
	
	C. Meyer (\cite{Meyer}) studied extensively double octics, in particular he gave a list of invariants of 450 examples of double octics.  Double octics with the Hodge number $h^{1,2}\le1$ where completely classified in \cite{C-KC}.

	\section{Elliptic fibrations determined by fourfold and fivefold points}
	Every point of an octic arrangement of multiplicity 4 or 5 determines an elliptic fibration of the double octic Calabi-Yau threefold.  Assume that the planes $P_{1},\dots,P_{4}$ from an octic arrangement $D=P_{1}\cup\dots\cup P_{8}$ intersect and denote $A=P_{1}\cap\dots\cap P_{4}$. Let $L$ be a general line in $\PP^{3}$ that contains the point $A$. Then the intersection $D.L$  equals $4A+Q_{1}+Q_{2}+Q_{3}+Q_{4}$, where $Q_{i}$ are distinct points. The double covering $\tilde E$ of $L$ branched along $L.D$ is singular curve of genus 1. Since in the process of resolution of a double cover of $\PP^{3}$ branched along $D$ we blow-up $\PP^{3}$ at $A$ and then replace the branch divisor by its pull-back minus four time exceptional divisor, the curve $\tilde E$ lifts to the double octic defined by $D$ as its normalization $E$ equal the double cover of $L$ branched along the divisor $L.D-4A = Q_{1}+\dots+Q_{4}$. In particular $E$ is an elliptic curve, and we get an elliptic fibration of the double octic.
	
	Following \cite{czarnecki} we can compute an equation of the elliptic fibration and its $j$-function. For simplicity, assume that $A=(0,0,0,1)$, we can write the equation of the plane $P_{i}$ in the form
	\[F_{i}(x,y,z)+G_{i}t\]
	where $F_{i}(x,y,z)$ is a polynomial linear with respect to variables $x,y,z$ and $G_{i}$ is a constant. Reordering. if necessary, planes $P_{i}$, we can assume that $G_{1}=\dots=G_{4}=0$ and $G_{5}=0$ when $A$ is a fivefold point and $G_{5}=1$ when $A$ is fourfold.

\begin{theorem}\label{thm:ellfive}
			If $A=(0:0:0:1)$ is a fivefold point then
	the double octic Calabi-Yau threefold defined by
	\begin{equation}\label{docsing5}
		u^{2}=F_{1}F_{2}F_{3}F_{4}F_{5}
		(F_{6}+t)(F_{7}+t)(F_{8}+t)
	\end{equation}
	is birational to the elliptic fibration
	\[u^{2} = \left(F_{1}F_{2}F_{3}F_{4}F_{5} F_{6}   +t \right) \left(F_{1}F_{2}F_{3}F_{4}F_{5} F_{7}  +t \right) \left(F_{1}F_{2}F_{3}F_{4}F_{5}F_{8}  +t \right).
	\]
	with discriminant
	\[
	\Delta = 16 (F_{1}F_{2}F_{3}F_{4}F_{5})^{6} \left({F_{7}}  -{F_{8}}  \right)^{2} 
	\left({F_{6}}  -{F_{8}}  \right)^{2} \left({F_{6}} -{F_{7}} \right)^{2}\]
	and the $J$-invariant 
	\[J= \frac{4 \left(
		{F_{6}}^{2} + {F_{7}}^{2}  + {F_{8}}^{2} 
		-{F_{6}} F_{7} 
		-{F_{6}} F_{8}  
		-{F_{7}} F_{8} 
		\right)^{3}}
	{27 {\left({F_{7}}  -{F_{8}} \right)^{2} \left({F_{6}} -{F_{8}} \right)^{2} \left({F_{6}}  -{F_{7}} \right)^{2}}}
	\]
	
\end{theorem}

	Let $a,b,c,d$ be four distinct complex numbers, 
the genus 1 curve 
\begin{equation}\label{eq:gen1}
	u^{2} = (t+a)(t+b)(t+c)(t+d)	
\end{equation}
become an elliptic curve after moving one of the points $(a,0)$, $(b,0)$, $(c,0)$ and $(d,0)$ to infinity. For instance the map 
\[(u,t)\longmapsto \left(\frac{(a-d)(b-d)(c-d)}{t+d}, \frac{(a-d)(b-d)(c-d)}{(t+d)^{2}}u\right)\] 
is a birational transformation of that curve onto the elliptic curve
\[u^{2} = (t+(b-d)(c-d))(t+(a-d)(c-d))(t+(a-d)(b-d)).\]
As a consequence, we get the following theorem

	\begin{theorem}\label{thm:ellfour}
				If $A=(0:0:0:1)$ is a fourfold point, then
		the double octic Calabi-Yau threefold defined by
		\begin{equation}\label{docsing4}
			u^{2}=F_{1}F_{2}F_{3}F_{4}(F_{5}+t)
			(F_{6}+t)(F_{7}+t)(F_{8}+t)
		\end{equation}
		is birational to the elliptic fibration
		
		\begin{multline*}
			u^{2}  = \Bigl(F_{1}F_{2}F_{3}F_{4} (F_{5} - F_{7}) (F_{5} - F_{6}) + t\Bigr) \times\\
			\Bigl(F_{1}F_{2}F_{3}F_{4} (F_{5} - F_{8}) (F_{5} - F_{6}) + t\Bigr)\Bigl(F_{1}F_{2}F_{3}F_{4} (F_{5} - F_{8}) (F_{5} - F_{7}) + t\Bigr)			
		\end{multline*}	
		with the discriminant
		\[\Delta = 16 (F_{1}F_{2}F_{3}F_{4})^{6} 
		\prod\limits_{5\le i<j\le 8}\left({F_{i}}  -{F_{j}}  \right)^{2} 
\]
		and the $J$-invariant
		
		\[
			J = \frac
			{4\left(6F_{5}F_{6}F_{7}F_{8}  -\sum\limits_{5\le i<j<k\le 8}F_{i}F_{j}F_{k}(F_{i}+F_{j}+F_{k}) + \sum\limits_{5\le i<j\le 8}F_{i}^{2}F_{j}^{2}
\right)^{3}}
{\left(\prod\limits_{5\le i<j\le 8}(F_{i}-F_{j})^{2}\right)}
\]
		
	\end{theorem}
	
		Elliptic fibrations resulting from the above theorem can be written as
		\[u^{2}=t^{3}+f_{6}(x,y,z)t^{2}+f_{12}(x,y,z)t+f_{18}(x,y,z),\]
		where $f_{6}, f_{12}, f_{18}$ are homogeneous polynomials in variables $x,y,z$ of degree 6, 12 and 18 respectively. Discriminant of the elliptic fibration is a non-reduced curve of degree 36, its reduction is a union of at most 10 lines (in the case of a fivefold point it is a union of at most 8 lines). Singular fibers at generic points of a line in the discriminant has Kodaira type that can be read from multiplicities of $\Delta$ and coefficients of the Weierstrass equation along that line.

	Elliptic fibration \eqref{docsing5}  has the same $J$-function as the fibration  	
	\[			u^{2}  = \Bigl( (F_{5} - F_{7}) (F_{5} - F_{6}) + t\Bigr) 
	\Bigl( (F_{5} - F_{8}) (F_{5} - F_{6}) + t\Bigr)\Bigl((F_{5} - F_{8}) (F_{5} - F_{7}) + t\Bigr)			
	\]
	but they  are not birational in general because they have singular fibers of different types at generic points of $F_{1}F_{2}F_{3}F_{4} = 0$.
	Smooth fibers of both fibrations are isomorphic, they differ by the quadratic twist by 
	$\sqrt{F_{1}F_{2}F_{3}F_{4}}$. 
	Similarly, for the elliptic fibrations ~\eqref{docsing5} and the elliptic fibration
	\[
		u^{2}=(F_{6}+t)(F_{7}+t)(F_{8}+t)		
	\]
	
%

\section{Elliptic fibrations determined by two skew double lines}

In this section we shall describe an elliptic fibration on a double octic Calabi-Yau threefold determined by four arrangement planes $P_{1}, P_{2}, P_{3}$ and $P_{4}$ that do not intersect at a point. The lines $L_{1,2}=P_{1}\cap P_{2}$ and $L_{3,4}=P_{3}\cap P_{4}$ are disjoint (skew lines in $\mathbb P^{3}$). Every point $Q\in \PP^{3}\setminus (L_{1,2}\cup L_{3,4})$ determine a unique line $L\subset \PP^{3}$ such that $L_{1,2}\cap L\not=\emptyset$ and $L_{3,4}\cap L\not=\emptyset$, in fact such a line is uniquely determined by the intersection points
$\{A\}:=L_{1,2}\cap L, \ \{B\}:=L_{3,4}\cap L$.

The intersection of the line $L$ with the branch divisor $D$ equals
\[L.D= 2A+2B+Q_{1}+Q_{2}+Q_{3}+Q_{4},\]
where $Q_{1}, Q_{2}, Q_{3}, Q_{4}$ are the intersection points of $L$ with the planes $P_{5},\dots,P_{8}$. 
The double covering $\tilde E$ of $L$ branched along the divisor $L.D$ is a singular curve of genus 1, since in the resolution process we blow-up $\PP^{3}$ at lines $L_{1,2}$ and $L_{3,4}$ the curve $\tilde E$ lifts to the double octic defined by $D$ as its normalization $E$ which is the double cover of $L$ branched along 
$L.D-2A-2B = Q_{1}+\dots+Q_{4}$. The curve $E$ is an elliptic curve and we get an elliptic fibration of a double octic.

To simplify further computations assume that $P_{1}=x, P_{2}=y, P_{3}=z, P_{4}=t$. The line $L$ that intersect  lines $L_{1,2}=\{x=y=0\}$ and $L_{3,4}=\{z=t=0\}$ at points $A=(0,0,\alpha_{0},\alpha_{1})$ and $B=(\beta_{0},\beta_{1},0,0)$ is given by the equations
\[L:\quad \alpha_{1}z-\alpha_{0}t = 0,\ \beta_{1}x-\beta_{0}y=0\]
so the curve $E$ is isomorphic to a curve given by the following equation
\[u^{2} = \frac1{(\alpha_{1}\beta_{0})^{4}} \prod_{i=5}^{8																																																																																										}\left(F_{i}(\alpha_{1}\beta_{0},\alpha_{1}\beta_{1},0,0)x + F_{i}(0,0,\alpha_{0}\beta_{0},\alpha_{1}\beta_{0})t\right)
\]
The above elliptic fibration is defined over the base of $L_{1,2}\times L_{3,4}$, if we prefer a fibration over $\PP^{2}$, then we dehomogenize it by letting $\alpha_{1}=\beta_{0}=1$ and then homogenize it back
obtaining

\[u^{2} = \prod_{i=5}^{8} \left(F_{i}(\gamma,\beta,0,0)x + F_{i}(0,0,\alpha,\gamma)t\right).
\]
By the substitution $\left(\frac 1t-\frac{F_{5}(\gamma,\beta,0,0)} {F_{5}(0,0,\alpha,\gamma)},\frac{u}{t^{2}} \right)$ in a similar manner as in the case of equation \eqref{eq:gen1}, and then using admissible transformation we can transform the equation into much simpler form, to improve legibility of appearing formulas we use the following notation
\[F_{i}:=F_{i}(\gamma,\beta,0,0), \quad G_{i}:=F_{i}(0,0,\alpha,\gamma).\]
\begin{theorem}
A double octic Calabi-Yau threefold defined by the octic arrangement
\[D:=\{(x,y,z,t)\in\PP^{3}: xyzt\cdot F_{5}(x,y,z,t)\cdots F_{8}(x,y,z,t)=0\}\]
is birational to the elliptic fibration
\begin{equation}\label{eq:gen2}
u^{2}=t\biggl(t-(F_{6} G_{8} - F_{8} G_{6}) (F_{5} G_{7} - F_{7} G_{5})\biggr) \biggl(t-(F_{6} G_{7} - F_{7} G_{6}) 
(F_{5} G_{8} - F_{8} G_{5})\biggr),	
\end{equation}
with the discriminant 
\[
\Delta=16\prod_{5\le i<j\le8}(F_{i}G_{j}-F_{j}G_{i})^{2} 
\]
and the $J$-invariant
\begin{multline*}
	J=\dfrac{\biggl(4\sum\limits_{\{i,j,k,l\}=\{5,6,7,8\},\atop i<j,\; k<l}F_{i}^{2}F_{j}^{2}G_{k}^{2}G_{l}^{2}
	-\sum\limits_{\{i,j,k,l\}=\{5,6,7,8\},\atop j<k} F_{i}^{2}F_{j}F_{k}G_{j}G_{k}G_{l}^{2}
	+6 {F_5} {F_6} {F_7} {F_8} {G_5} {G_6} {G_7} {G_8}
\biggr)^{3}}
	{\Biggl(27 \prod\limits_{5\le i<j\le8}\bigl(F_{i}G_{j}-F_{j}G_{i}\bigr)^{2}\Biggr)}.
\end{multline*}
\end{theorem}

%
%
%
%
\section{K3 fibrations on a double octic}
In this section we shall describe three types of K3 fibrations on a double octic:  Kummer fibration determined by two disjoint, intersecting quadruples of arrangement planes, K3 fibrations determined by a pencil of planes containing an arrangement double or triple line, K3 fibration determined by two disjoint arrangement double or triple lines.

\subsection{Kummer fibrations}

Assume that the arrangement contains two disjoint, intersecting quadruples, after renumbering we can assume that 
\[P_{1}\cap\cdots\cap P_{4} = \{Q_{1}\} \quad \text{ and } P_{5}\cap\cdots\cap P_{8} = \{Q_{2}\},\]
($Q_{1}$ and $Q_{2}$ are arrangement fourfold or fivefold points).
Let $m$ be the pencil (a projective line) of planes containing points $Q_{1}$ and $Q_{2}$. Take a generic plane $S\in m$ and let $m_{1}$  be the projective line (pencil) of lines in $S$ containing point $Q_{1}$. The double cover of $m_{1}$ branched along $P_{1}\cap S, \dots, P_{4}\cap S$ is an elliptic curve. 
We get an elliptic surface $S_{1}\lra m$. 

Considering the planes $P_{5},\dots,P_{8}$ and their intersection point $Q_{2}$ we construct in a similar manner an elliptic surface $S_{2}\lra m$.

Surfaces $S_{1}$ and $S_{2}$ are clearly birational to a double covering 
of a projective plane branched along a union of four lines. If we denote by $\Pi_{i}$ the projective plane of lines containing $Q_{1}$, then the planes $P_{1}, \dots, P_{4}$ define four lines in $\Pi_{1}$, the surface $S_{1}$ is birational to the double cover of $\Pi_{1}$ branched along these four lines, we have a similar description for the surface $S_{2}$. 

If we assume that $Q_{1}=(0:0:0:1)$ and $Q_{2}=(1:0:0:0)$, then planes $P_{1},\dots,P_{4}$ have equations of the form $F_{1}(x,y,z),\dots,F_{4}(x,y,z)$. planes $P_{5},\dots,P_{8}$ have equations of the form $F_{5}(y,z,t),\dots,F_{8}(y,z,t)$ and surfaces $S_{1}$ and $S_{2}$ are birational to respectively
\[v_{1}^{2} = F_{1}(x,y,z)\cdots F_{4}(x,y,z) \quad \text{ and } \quad 
v_{2}^{2} = F_{5}(y,z,t)\cdots F_{8}(y,z,t).\]
The (fiber-wise) elliptic involutions on $S_{1}$ and $S_{2}$ are given by
\[i_{1}:S_{1}\ni(x,y,z,v_{1})\longmapsto (x,y,z,-v_{1})\in S_{1}
\]
and \[i_{2}:S_{2}\ni(y,z,t,v_{2})\longmapsto (y,z,t,-v_{2})\in S_{2}.\]
The double octic is birational to the quotient
\[(S_{1}\times_{\PP^{1}}S_{2}) / (i_{1}\times i_{2})\]
of the fiber product of $S_{1}$ and $S_{2}$ by the product of $i_{1}$ and $i_{2}$.  

\begin{theorem}
	The double octic Calabi-Yau threefold 
	\[u^{2}=F_{1}(x,y,z)\cdots F_{4}(x,y,z) F_{5}(y,z,t)\cdots F_{8}(y,z,t)\] 
	is birational to the Kummer fibraion for elliptic surfaces 
	\[v_{1}^{2}=F_{1}(x,y,z)\cdots F_{4}(x,y,z) \quad \text{ and } \quad v_{2}^{2}F_{5}(y,z,t)\cdots F_{8}(y,z,t).\]
\end{theorem}
Positions and types of singular fibers can be easily computed, there are six types of elliptic surfaces appearing in this construction (presented by eight types of double quartic surfaces), for details see \cite{C-KC}. 

More than half of families of double octics admits a Kummer fibration, considering the Kummer fibrations we found a large number of birational maps or more generally correspondences  between different families of double octics.
The $J$-function and discriminant for elliptic surfaces $S_{1}$ and $S_{2}$ can be easily computed as in Theorem~\ref{thm:ellfour} and Theorem~\ref{thm:ellfive}.

\subsection{K3 fibrations by double sextics}.

For any pair $P_{i}, P_{j}$, $i\not=j$ of arrangement planes the planes through the double or triple line $L=P_{i}\cap P_{j}$ determine a K3 fibration of the double octic. Namely, for a plane $P_{s_{0}:s_{1}}$ defined by the polynomial $s_{0}F_{i}+s_{1}F_{j}$, where $F_{i}$ and $F_{j}$ are equations of planes $P_{i}$ and $P_{j}$ respectively, the fiber is given by the double cover of $P_{s_{0}:s_{1}}$ branched along a union of six lines
\[l_{k}:=P_{(s_{0}:s_{1})}\cap P_{k}, \quad k\not\in\{i,j\}.\]
Double cover of the projective plane branched along a smooth curve of degree 6 is a K3 surface, if the branch curve is singular (possibly reducible) but has at worst Kleinian singularities, then the double cover has also Kleinian singularities (of the same type) and its minimal resolution of singularities is a K3 surface. 
This assumption is satisfied in particular, when the branch curve is a union of six lines, with at most three intersecting at a point (for details see  \cite{Persson}), 
Since for a generic choice of $(s_{0}:s_{1})\in\PP^{1}$ the lines $l_{k}$ are all distinct and no four of them intersect, the double cover has a resolution of singularities which is a K3 surface.

To get an explicit equation of the K3 fibration assume that $F_{7}=z$, $F_{8}=t$ and the line $L$ is given by $z=t=0$. For any $s\in\CC$ the plane $P_{-s:1}$ is given by the equation $t=sz$

\begin{theorem}
	A double octic Calabi-Yau threefold defined by the octic arrangement
	\[D:=\{(x,y,z,t)\in\PP^{3}:P_{1}(x,y,z,t)\cdots P_{6}(x,y,z,t)zt=0\}\]
	is birational to the K3 fibration 
	\[u^{2}=F_{1}(x,y,z,sz) \cdots F_{6}(x,y,z,sz)\]
\end{theorem}
The fiber at infinity is the double sextic 
\[u^{2}=F_{1}(x,y,0,t)\cdots F_{6}(x,y,0,t).\]

Arrangement double and triple lines disjoint from $L$ yield double and triple points of the 
branch curve in a generic fiber, similarly - fourfold and fivefold points lying on $L$. The special fibers depend on double and triple lines intersecting $L$ and fourfold and fivefold points not lying on $L$. 
Explicit description of general and singular fibers was given by P. Borówka, who studied in details the case of rigid double octics  (\cite{Bor}).

\subsection{K3 fibrations by double quadrics}.

A pair of skew arrangement double or triple lines determines a K3 fibration of a double octic. The lines $L_{ij}=P_{i}\cap P_{j}$ and $L_{kl}=P_{k}\cap P_{l}$, $\#\{i,j,k,l\}=4$, are skew exactly when $P_{i}\cap P_{j}\cap P_{k}\cap P_{l}=\emptyset$. The base locus of the quadric pencil $Q_{\lambda_{0}:\lambda_{1}}$, $(\lambda_{0}:\lambda_{1})\in\PP^{1}$ defined by quadratic forms $\lambda_{0}F_{i}F_{k}-\lambda_{1}F_{j}F_{l}$
is the union of four lines $L_{ij}$, $L_{il}$, $L_{kj}$ and $L_{kl}$. In the resolution process we blow-up these lines, consequently we get a fibration on the double octic. Generic fiber $S_{\lambda_{0}:\lambda_{1}}$ is a resolution of the double cover of the quadric $G_{\lambda_{0}:\lambda_{1}}$ branched along the intersection of $Q_{\lambda_{0}:\lambda_{1}}$ with the octic arrangement.
The intersection equals 
\[S_{\lambda_{0}:\lambda_{1}}D = 2L_{ij}+2L_{il}+2L_{kj}+2L_{kl}+\sum_{m\not\in\{i,j,k,l\}} Q_{\lambda_{0}:\lambda_{1}}\cap P_{m}.\]
Consequently, the branch locus 
\[
S_{\lambda_{0}:\lambda_{1}}D - (2L_{ij}+2L_{il}+2L_{kj}+2L_{kl}) = \sum_{m\not\in\{i,j,k,l\}} Q_{\lambda_{0}:\lambda_{1}}\cap P_{m}.
\]
is a curve of type $(4,4)$, while the canonical divisor of a quadric has type $(-2,-2)$. Consequently, the resolution of the double cover of $Q_{\lambda_{0}:\lambda_{1}}$ is a K3 surface. 

The quadric $Q_{\lambda_{0}:\lambda_{1}}$ is isomorphic to $\PP^{1}\times \PP^{1}$, an explicit isomorphism can be given by letting

\[\lambda_{0} p_{0}q_{0}=F_{i}, \ \ \lambda_{1} p_{1}q_{1}=F_{k}, \ \ \lambda_{0}p_{0}q_{1}=F_{j}, \ \ \lambda_{0}p_{1}q_{0}=F_{l}.\]
The assumption that lines $L_{ij}$ and $L_{kl}$ are disjoint implies that the linear map $(F_{i},F_{j},F_{k},F_{l})$ is an isomorphism. Consequently, the above system of equations have a unique solution in $x,y,z,t$, substituting the solution into double octic, yield a family of double quadrics.
To simplify computations we assume that the double lines equal  $x=y=0$ and $z=t=0$.

\begin{theorem}
	The double octic
	\[u^{2} = xyztF_{5}(x,y,z,t)\dots F_{8}(x,y,z,t)\]
	is birational to a double quadrics fibration
	\[u^{2} = \lambda F_{5}(p_{0}q_{0}, p_{0}q_{1}, p_{1}q_{0}, \lambda p_{1}q_{1}) \cdots F_{8}(p_{0}q_{0}, p_{0}q_{1}, p_{1}q_{0}, \lambda p_{1}q_{1}).\] 
\end{theorem}

Fibrations by double quadrics are  more difficult to classify then fibrations by Kummer surfaces or double sextics. A basic invariant of a family of K3 surfaces is given by the Picard-Fuchs operator, which is a differential operator annihilating the period integral (\cite{CvS4}). Unfortunately, the Picard-Fuchs operators of families of K3 surfaces do not have as good properties as in the case of Calabi-Yau threefolds. In particular, we do not have a nice classification of the possible Riemann Symbols. Nevertheless a birational isomorphism of two fibrations results in equal Picard-Fuchs operators, the opposite implication also is expected to hold true. 

A singular point of the branch divisor determines an elliptic fibration on the double quadric, checking a fibration preserving isomorphism of two elliptically fibered K3 surfaces is an elementary linear algebra. 
In the next section we shall use the above simple observation to find a birational transformation between two families of double octic.

%
%
%
	\section{Application}
	Consider two one parameter families of double octics defined by the following equations
	\begin{multline*}
		\text{\bf Arrangement No. 35:}\\
	X^{35}_{A,B}:\quad u^{2}=(Ax + By)t(Ax + By + At)xy(By + Bz + At)(x + y + z + t)z	
	\end{multline*} 
	and
	\begin{multline*}
		\text{\bf Arrangement No. 71:}\\
	X^{71}_{A,B}: \quad u^{2}=yx(x + y)(By-Az-At)(x + y + z + t)z(Ax-By + Az)t
	\end{multline*} 
	These families have the following Picard-Fuchs operators

	\textbf{35:\quad }
	\begin{center}
		\mbox{\(\displaystyle\Theta^{2}( \Theta -\tfrac12)( \Theta +\tfrac12)\)}
		\mbox{\(\displaystyle-t( \Theta +\tfrac12)^{2}(2 \Theta^{2}+2 \Theta +1)\)}
\mbox{\(\displaystyle+t^{2}(\Theta +1)^{2}(5 \Theta^{2}+10 \Theta +6)\)}

\mbox{\(\displaystyle-2t^{3}(\Theta +1)^{2}(\Theta +2)^{2}\)}
	\end{center}

	with Riemann Symbols
	
	\(\rule[1.7cm]{56mm}{0mm}\left\{ \begin {array}{cccc}0&1/2&1&\infty  \\
	\hline -1/2&0&0&1\\
	0&1&1/2&1\\
	0&1&1/2&2\\
	1&2&1&2
	\end {array} \right\} \)
	
	\medskip
	
	\textbf{71:\quad}
	\(
	\Theta\, \left( \Theta-1 \right)  \left( \Theta -\frac12 \right) ^{2}
	+t{\Theta}^{2} \left( 4\,{\Theta}^{2}+1 \right) 
	+\frac1{16}\,{t}^{2} \left( 20\,{\Theta}^{2}+20\,\Theta+9 \right)  \left( 2\,\Theta+1 \right) ^{2}
		\)

	\(\rule[6mm]{4cm}{0cm}	
	+\frac18\,{t}^{3} \left( 2\,\Theta+3 \right) ^{2} \left( 2\,\Theta+1 \right) ^{2}
	\)

with Riemann Symbols

	\(\rule[1.7cm]{56mm}{0mm}\left\{
	\begin {array}{cccc}
	0&-1/2&-1&\infty\\\hline
	0&0&0&1/2\\
	1/2&1&1/2&1/2\\
	1/2&1&1/2&3/2\\
	1&2&1&3/2
	\end {array} \right\} \)

The two differential operators coincide up-to the variable sign $\phi:t\longmapsto-t$ and multiplication with an algebraic function, moreover elements of both families have the same Hodge numbers $h^{1,1}=49$, $h^{1,2}=1$ and hence also the topological Euler characteristic $\chi=-96$. Consequently, the two families are expected to be birational (\cite{CvS3}). 
In this paper we shall work with alternative equations for these families (obtained from a more systematic studies of double octics)
	\begin{eqnarray*}
	X_{A_{0}:A_{1}}:&&
	u^{2}=(x - y)xy(y + t)(x + t)z(A_{0}x + A_{1}y + (-A_{0} + A_{1})z + A_{1}t)(z + t)\\
	Y_{A_{0}:A_{1}}:&&
	u^{2}=xy(x + y)(A_{0}x + A_{1}y + A_{1}z)zt((-A_{0} + A_{1})x + A_{1}y + A_{1}t)(x + y + z + t)
\end{eqnarray*} 
and related by the following changes of coordinates
\begin{eqnarray*}
	\mathcal X^{35}\ni (x,y,z,t,u,A,B)&\longmapsto&(-By,Ax,Bz,By+At,(A^{3}B^{3})^{1/2}u,-A+B,B)
	\in\mathcal X
	\\
	\mathcal X^{71}\ni(x,y,z,t,u,A,B)&\longmapsto&(y,x,-(x+y+z+t),z,u,A+B,A)\in\mathcal Y
\end{eqnarray*}
We shall use K3 fibrations in double quadrics described in the last section to derive an explicit birational map between elements of families $\mathcal{X}$ and $\mathcal{Y}$ 
\begin{theorem}
	The following map is birational 
	\begin{multline*}
		\mathcal{X}\ni(x,y,z,t,u,A_{0},A_{1}) \longmapsto
		\Bigl(
		-(y-z)(x-y)x(A_{0}-A_{1}),\\ -(x-y)((A_{0}-A_{1})xz-(A_{0}-A_{1})yz+A_{1}ty+A_{1}y^2),\\ (y-z)(x-y)(A_{0}x-(A_{0}-A_{1})z+A_{1}t+A_{1}y),\\ (x+t)(A_{0}xz-A_{0}yz+A_{1}ty-A_{1}xy+A_{1}
		y^2+A_{1}yz),\\ (A_{0}-A_{1})\left(\frac{A_{0}-A_{1}}{A_{0}}\right)^{1/2}
		(A_{0}x-A_{0}z+A_{1}t-A_{1}x+A_{1}y+A_{1}z)\\
		(x-y)^2(y-z)(A_{0}xz-A_{0}yz+A_{1}ty-A_{1}xz+A_{1}y^2+A_{1}yz)\\
		(A_{0}xz-A_{0}yz+A_{1}ty-A_{1}xy+A_{1}y^2+A_{1}yz), A_{0}, A_{0}-A_{1}
		\Bigr)\in\mathcal{Y}
	\end{multline*}
	The inverse birational map is given by
	\begin{multline*}
		(x,y,z,t,u,A_{0},A_{1}) \longmapsto\\
		\Bigl(
		x(A_{0}^2xy-A_{0}A_{1}tx-A_{0}A_{1}ty-A_{0}A_{1}x^2-3A_{0}A_{1}xy-A_{0}A_{1}xz-A_{0}A_{1}y^2\\+A_{1}^2tx+A_{1}^2ty+A_{1}^2x^2
		+3A_{1}^2xy+A_{1}^2xz+2A_{1}^2y^2+A_{1}^2yz),\\
		-(x+y)A_{1}(A_{0}tx+A_{0}x^2+A_{0}xz-A_{1}tx-A_{1}x^2-A_{1}xy-A_{1}xz-A_{1}yz),\\ -(A_{0}x+A_{1}y+A_{1}z)(A_{0}tx+A_{0}x^2+A_{0}xz-A_{1}tx-A_{1}x^2-A_{1}xy-A_{1}xz-A_{1}yz),\\
		-A_{0}^2x^2y+A_{0}A_{1}tx^2+2A_{0}A_{1}txy+A_{0}A_{1}x^3+3A_{0}A_{1}x^2y+A_{0}A_{1}x^2z+A_{0}A_{1}xy^2\\-A_{1}^2tx^2-2A_{1}^2txy+A_{1}^2tyz-A_{1}^2x^3-3A_{1}^2x^2y-A_{1}^2x^2z-2A_{1}^2xy^2-A_{1}^2xyz, A_{0}, A_{0}-A_{1},\\
		A_{1}\sqrt{A_{0}A_{1}}y(A_{0}x-A_{1}x+A_{1}z)^2
		(A_{0}tx+A_{0}x^2+A_{0}xz-A_{1}tx-A_{1}x^2-A_{1}xy-A_{1}xz-A_{1}yz)\\
		(A_{0}^2xy-A_{0}A_{1}tx-A_{0}A_{1}ty-A_{0}A_{1}x^2-3A_{0}A_{1}xy-A_{0}A_{1}xz-A_{0}A_{1}y^2+A_{1}^2tx+A_{1}^2ty+A_{1}^2x^2\\+3A_{1}^2xy+A_{1}^2x
		z+2A_{1}^2y^2+A_{1}^2yz)(A_{0}x-A_{1}x-A_{1}y)
		\Bigr)
	\end{multline*}
\end{theorem}

The prove of this theorem could be accomplished by substituting the equations of the above map into equation of family $\mathcal{Y}$.
We shall instead describe in full details how the above formulas were derived .

For both families we consider all pairs of skew (disjoint) arrangements double lines and compute the corresponding families of double covers of $\PP^{1}\times \PP^{1}$ branched along a curve of type $(4,4)$. This way we get 165 K3 fibrations for generic elements of both families. As the expected birational transformation should map $X_{2:1}$ to $Y_{2:1}$, we compute the Picard-Fuchs operators for all constructed K3 fibrations on both double octics. 

The choice of double lines $y=y+t=0$ and $x+t=z=0$ for the family No. 35 yields the following family of double quadrics 
\begin{multline*}
	\tilde {\mathcal X}: \quad v^{2}=lq_{1}\left(l p_{1} q_{1} - p_{0} q_{0} - p_{0} q_{1}\right)\left(l p_{1} q_{1} - p_{0} q_{0} - p_{1} q_{0}\right) \\
	\left(A_{0} p_{0} q_{1} - A_{0} p_{1} q_{0} + A_{1} p_{0} q_{0} + A_{1} p_{1} q_{0}\right) 
	\left(l p_{1} - p_{0}\right) 
\end{multline*}
and the birational map 
\begin{multline*}
((p_{0}:p_{1}),(q_{0}:q_{1}), l,v)\longmapsto \\
\left(p_{0}q_{1}+p_{0}q_{0}-lp_{1}q_{1}, 
p_{0}q_{0}, 
p_{1}q_{0}+ p_{0}q_{0}-lp_{1}q_{1}, -p_{0}q_{0}+lp_{1}q_{1},
A_{0}^{(7/2)}\sqrt{-1}p_{0}p_{1}q_{0}q_{1}v\right)	
\end{multline*}
to the double octic $X_{A_{0}:A_{1}}$.
The Picard-Fuchs operator of the K3 fibration on the fiber at $A_{0}=2, A_{1}=1$, computed with the package of P. Lairez \cite{Lai},
equals  
\begin{center}
	\mbox{\(\displaystyle\frac{1}{2}\Theta^{3}(\Theta -1)(2 \Theta -1)\)}
	\mbox{\(\displaystyle+\frac{1}{2^{2}}t\Theta(64 \Theta^{4}-16 \Theta^{3}+57 \Theta^{2}+13 \Theta +2)\)}
	\mbox{\(\displaystyle+\frac{1}{2^{3}}t^{2}(1626 \Theta^{5}-1929 \Theta^{4}+2971 \Theta^{3}+999 \Theta^{2}+434 \Theta +84)\)}
	\mbox{\(\displaystyle-\frac{1}{2^{4}}t^{3}(30878 \Theta^{5}-43461 \Theta^{4}+38565 \Theta^{3}-1644 \Theta^{2}+325 \Theta +414)\)}
	\mbox{\(\displaystyle-\frac{3}{2^{5}}t^{4}(117164 \Theta^{5}-87386 \Theta^{4}+59195 \Theta^{3}-53267 \Theta^{2}-33831 \Theta -8666)\)}
	\mbox{\(\displaystyle+\frac{1}{2^{6}}t^{5}(2462562 \Theta^{5}-272259 \Theta^{4}+2169118 \Theta^{3}-1517127 \Theta^{2}-1389211 \Theta -421374)\)}
	\mbox{\(\displaystyle-\frac{1}{2^{7}}t^{6}(11597288 \Theta^{5}+1250850 \Theta^{4}+13853841 \Theta^{3}-4331658 \Theta^{2}-6722507 \Theta -2503806)\)}
	\mbox{\(\displaystyle+\frac{3}{2^{7}}t^{7}(6448356 \Theta^{5}+397704 \Theta^{4}+2247035 \Theta^{3}-9265571 \Theta^{2}-8508114 \Theta -2647368)\)}
	\mbox{\(\displaystyle-\frac{1}{2^{5}}t^{8}(5585898 \Theta^{5}-661155 \Theta^{4}-8099216 \Theta^{3}-24512394 \Theta^{2}-20618395 \Theta -6251490)\)}
	\mbox{\(\displaystyle+\frac{1}{2^{5}}t^{9}(4027956 \Theta^{5}-4534002 \Theta^{4}-23855270 \Theta^{3}-46360215 \Theta^{2}-36321487 \Theta -10566366)\)}
	\mbox{\(\displaystyle-\frac{1}{2^{5}}t^{10}(1311320 \Theta^{5}-9622850 \Theta^{4}-39222885 \Theta^{3}-64789870 \Theta^{2}-49012229 \Theta -14025450)\)}
	\mbox{\(\displaystyle+\frac{3^{2}}{2^{5}}t^{11}(\Theta +1)(47564 \Theta^{4}+1006284 \Theta^{3}+2876321 \Theta^{2}+3230796 \Theta +1246772)\)}
	\mbox{\(\displaystyle-\frac{3^{3}}{2^{3}}t^{12}(\Theta +1)(\Theta +2)(6246 \Theta^{3}+33141 \Theta^{2}+47833 \Theta +20649)\)}\\
	\mbox{\(\displaystyle+3^{3}17^{2}t^{13}(\Theta +1)^{3}(\Theta +2)(\Theta +3)\)}
\end{center}

In a similar manner the choice of double lines
$x+y =  (-A_{0} + A_{1}) x + A_{1} y + A_{1} t = 0$ and
$A_{0} x + A_{1} y + A_{1} z = t = 0$ for family No. 71 yields the following family of double quadrics
\begin{multline*}
\tilde{\mathcal Y}:\quad	v^{2} =  l(A_{0}lp_{1}q_{1}+A_{0}p_{0}q_{0}-A_{1}p_{0}q_{0}-A_{1}p_{1}q_{0})(A_{0}lp_{1}q_{1}-A_{1}p_{0}q_{0}-A_{1}p_{1}q_{0})\\(A_{0}^2lp_{1}q_{1}-A_{0}A_{1}lp_{1}q_{1}+A_{0}^2p_{0}q_{1}-2A_{0}A_{1}p_{0}q_{0}-A_{0}A_{1}p_{1}q_{0}+A_{1}^2p_{0}q_{0}+A_{1}^2p_{1}q_{0})\\(A_{0}^2lp_{1}q_{1}-A_{0}A_{1}lp_{1}q_{1}+A_{0}^2
p_{0}q_{1}-A_{0}A_{1}p_{0}q_{0}+A_{1}^2p_{0}q_{0}+A_{1}^2p_{1}q_{0})
\end{multline*}
and the birational map
\begin{multline*}
(p_{0},p_{1}, q_{0},q_{1}, l,v)\longmapsto 
\Bigl(\frac1{{A_{0}}}\left({A_{1}} {p_{0}} {q_{0}} +{A_{1}} {p_{1}} {q_{0}} -A_{0}l {p_{1}} {q_{1}}\right), 
\frac1{A_{0}}(A_{0}lp_{1}q_{1}+(A_{0}-A_{1})p_{0}q_{0} - A_{1}p_{1}q_{0}),\\ \frac1{A_{0}A_{1}}\bigl(A_{0}(A_{0}-A_{1})lp_{1}q_{1}+A_{0}^2p_{0}q_{1}-(2A_{0}-A_{1})A_{1}p_{0}q_{0}-(A_{0}-A_{1})A_{1}p_{1}q_{0}\bigr), p_{1}q_{0}, \sqrt{-1}\frac{A_{0}^{2}}{A_{1}}p_{0}p_{1}q_{0}q_{1}v \Bigr)
\end{multline*}

The Picard-Fuchs operator of the K3 fibration on the fiber at $A_{0}=2, A_{1}=1$
is identical as for the family No. 35. Consequenctly we search for a birational map $\tilde{\mathcal X}\lra \tilde{\mathcal Y}$ compatible with K3-fibrations.

Double quadrics $\tilde{\mathcal X}$ and $\tilde{\mathcal Y}$ are K3 surfaces, every double point of the branch locus determines an elliptic fibrations, for instance the point $p_{0}=q_{1}=0$ determines elliptic fibration 
\begin{multline*}
	w^{2} = p_{1}l(lp_{0}t_{0}-p_{0}t_{1}+p_{1}t_{1})(lp_{1}t_{0}-p_{0}t_{0}-p_{1}t_{1})\\
	(A_{0}lp_{0}p_{1}t_{0}-A_{1}lp_{0}p_{1}t_{0}-A_{0}p_{0}^2t_{0}-A_{0}p_{0}p_{1}t_{1}+A_{0}p_{1}^2t_{1}-A_{1}p_{1}^2t_{1})(lp_{1}-p_{0})t_{0}
\end{multline*}
from the pencil 
\[p_{0}q_{0}t_{0}=p_{1}q_{1}t_{1}.\]

In a similar manner for the family 71 fixing the double point $p_{1}=q_{0}=0$ we get an elliptic fibration
\begin{multline*}
	w^{2}=l(A_{_{0}}lq_{_{1}}t_{0}+A_{0}q_{1}t_{1}-A_{1}q_{0}t_{0}-A_{1}q_{1}t_{1})(A_{0}lq_{1}t_{0}-A_{1}q_{0}t_{0}-A_{1}q_{1}t_{1})\\	(A_{0}^2lq_{0}q_{1}t_{0}-A_{0}A_{1}lq_{0}q_{1}t_{0}+A_{0}^2q_{1}^2t_{1}-A_{0}A_{1}q_{0}^2t_{0}-2A_{0}A_{1}q_{0}q_{1}t_{1}+A_{1}^2q_{0}^2t_{0}+A_{1}^2q_{0}q_{1}t_{1})\\
	(A_{0}^2lq_{0}q_{1}
	t_{0}-A_{0}A_{1}lq_{0}q_{1}t_{0}+A_{0}^2q_{1}^2t_{1}-A_{0}A_{1}q_{0}q_{1}t_{1}+A_{1}^2q_{0}^2t_{0}+A_{1}^2q_{0}q_{1}t_{1})
\end{multline*}
from the pencil $p_{1}q_{1}t_{1}=p_{0}q_{0}t_{0}$.

Solving right-hand sides for $t_{0}$ we get
\begin{multline}\label{quadr35}
	0, \frac{p_{1}t_{1}}{lp_{1}-p_{0}}, \frac{t_{1}(p_{0}-p1)}{lp_{0}}, \frac{p_{1}t_{1}(A_{0}p_{0}-A_{0}p_{1}+A_{1}p_{1})} {p_{0}(A_{0}lp_{1}-A_{1}lp_{1}-A_{0}p_{0})}
\end{multline}
and 
\begin{multline}\label{quadr71}
		-\frac{q_{1}t_{1}(A_{0}-A_{1})}{A_{0}lq_{1}-A_{1}q_{0}}, \frac{A_{1}q_{1}t_{1}}{(A_{0}lq_{1}-A_{1}q_{0})}, -\frac{q_{1}t_{1}(A_{0}^2q_{1}-A_{0}A_{1}q_{0}+A_{1}^2q_{0})}
		{q_{0}(A_{0}^2lq_{1}-A_{0}A_{1}lq_{1}+A_{1}^2q_{0})},\\ -\frac{q_{1}t_{1}(A_{0}^2q_{1}-2A_{0}A_{1}q_{0}+A_{1}^2q_{0}
	)}{q_{0}(A_{0}^2lq_{1}-A_{0}A_{1}lq_{1}-A_{0}A_{1}q_{0}+A_{1}^2q_{0})}
\end{multline}
Cross-ratios of the above quadruples equal

\[\frac{(A_{0}p_{0}-A_{0}p_{1}+A_{1}p_{1})(lp_{1}^2+p_{0}^2-p_{0}p_{1})}{p_{0}^2(A_{1}lp_{1}+A_{0}p_{0}-A_{0}p_{1})}
\]
and
\[
\frac{(q_{0}-q_{1})(A_{0}^2lq_{1}^2-A_{0}A_{1}q_{0}q_{1}+A_{1}^2q_{0}^2)}
{q_{0}^2(A_{0}lq_{1}-A_{1}lq_{1}-A_{0}q_{1}+A_{1}q_{0})A_{1}}
\]
From the shapes of the above rational functions one easily deduces that the map
\[(p_{0},p_{1},q_{0},q_{1},A_{0},A_{1})\longmapsto (q_{0},q_{1},A_{0}p_{0},(A_{0}-A_{1})p_{1},A_{0}, A_{0}-A_{1})\]
transform the second cross-ratio into the first one. Consequently, this map transform the quadruple  ~\eqref{quadr71} into a quadruple projectively equivalent to \eqref{quadr35}, computing the M\"obius transformation between the two equivalent quadruples we get birational map between elliptic fibrations
\begin{multline*}
\bigl(p_{_{0}},p_{_{1}},q_{_{0}},q_{_{1}},A_{_{0}},A_{_{1}},t_{_{0}},t_{_{1}}\bigr)\longmapsto \biggl(q_{_{0}},q_{_{1}},A_{_{0}}p_{_{0}},(A_{_{0}}-A_{_{1}})p_{_{1}},A_{_{0}}, A_{_{0}}-A_{_{1}}, \\
\frac{(A_{_{0}}-A_{1})p_{1}(A_{0}lp_{1}^2t_{0}-2A_{1}lp_{0}p_{1}t_{0}-A_{0}p_{0}p_{1}t_{0}+A_{1}p_{0}^2t_{0}+A_{1}p_{0}p_{1}t_{1}-A_{1}p_{1}^2t_{1})t_{1}
}{A_{0}(A_{0}lp_{0}p_{1}t_{0}-2A_{1}lp_{0}p_{1}t_{0}-A_{0}p_{0}^2t_{0}-A_{0}p_{0}p_{1}t_{1}+A_{0}p_{1}^2t_{1}+A_{1}p_{0}^2t_{0}+A_{1}p_{0}p_{1}t_{1}-A_{1}p_{1}^2t_{1})(lp_{1}-p_{0})
},t_{1}\biggr)	
\end{multline*}

Composing computed maps we get the birational map  $\mathcal X\lra \mathcal Y$ in the assertion of the theorem.
In a similar way (exchanging the role of $\mathcal X^{35}$ and $\mathcal X^{71}$) we compute also the inverse birational  map.

\end{document}